\documentclass[12pt,oneside]{article}

\usepackage{amsfonts,amsmath,amssymb,amsthm,color,mathrsfs,nicefrac,stmaryrd,etexcmds,xfrac}
\usepackage[numbers]{natbib}
\usepackage{graphicx}
\usepackage[top=2in,bottom=1.5in,right=1in,left=1in]{geometry}
\usepackage[retainorgcmds]{IEEEtrantools}
\newcommand{\comment}[1]{}  

\newcommand{\Mrank}[1]{\mathrm{rank}\!\left(\right)}
\newcommand{\norm}[1]{\left\Vert#1\right\Vert}																				

\newcommand{\eqd}{\overset{\mathrm{d.}}{=}}																					

\newcommand{\real}{\mathbb{R}}

\DeclareMathOperator*{\argmin}{\arg\!\min}

\def\mby{\mathbf{y}}
\def\mbx{\mathbf{x}}
\def\wX{\widetilde{X}}
\def\wY{\widetilde{Y}}
\def\half{\hbox{$1\over2$}}

\newtheorem{thm}{Theorem}[section]

\newtheorem{cor}{Corollary}[section]
\newtheorem{prop}{Proposition}[section]

\begin{document}


\author{Novin Ghaffari \& S. G. Walker\footnote{
Department of Statistics \& Data Sciences, University of Texas at Austin,
and Department of Mathematics, University of Texas at Austin. s.g.walker@math.utexas.edu}}
\title{Parseval's Identity and Optimal Transport Maps\footnote{A correction to ``On Multivariate Optimal Transportation''}}	
	\date{}
	
	\maketitle
	
	\begin{abstract}
	Recent findings for optimal transport maps between distribution functions sharing the same copula show that componentwise the solution is the optimal map between marginal distributions. This is an important discovery since in the multivariate setting optimal maps are difficult to find and only known in a few special cases. In this paper, we extend the result on common copulas by showing that orthonormal transformations of variables sharing a common copula also have a known optimal map. We illustrate this by establishing optimal maps between members of a class of scale mixture of normal distributions.
	\end{abstract}
	
	
\vspace{0.2in}
\noindent	
{\it keywords}:
Positive definite matrix; Orthonormal matrix; Copula; Optimal transportation; Convex function; Barycenter.
	

	
\section{Introduction}
The Wasserstein metric defines a distance between probability measures and can be thought of as the minimal effort required to transport one probability measure to another under a cost function. In probability and statistics, the Wasserstein distance is typically studied under $L^p$ costs; $c(\mathbf{x},\mathbf{y})=\norm{\mathbf{x}-\mathbf{y}}_p$. The most common case is the quadratic or $2$-Wasserstein, with cost $\norm{\mathbf{x}-\mathbf{y}}_2^2$ and we will use this exclusively throughout the paper. 

The $W_2$ distance between probability measures $\mathbb{P}$ and $\mathbb{Q}$ is defined as  
	\begin{IEEEeqnarray}{RCL}\label{Weqn}
		W_2^2(\mathbb{P},\mathbb{Q})&=&\inf_{\pi\in \Pi(\mathbb{P},\mathbb{Q})}\int_{\real^{2d}} \norm{\mathbf{x}-\mathbf{y}}_2^2\,d\,\pi(\mathbf{x},\mathbf{y}),
	\end{IEEEeqnarray}
	where $\Pi(\mathbb{P},\mathbb{Q})$ denotes the set of probability distributions on $\real^{2d}$ with marginal probability distributions $\mathbb{P}$ and $\mathbb{Q}$. 
See \cite{villani2003} and  \cite{villani2006} for more details.
	
Given random variables $X\sim\mathbb{P}$ and $Y\sim\mathbb{Q}$ such that their joint distribution $(X,Y)\sim\pi^*$ attains the infimum of (\ref{Weqn}), we say that $X$ and $Y$ are \emph{optimally coupled} with respect to $W_2$. When the context is clear we simply refer to $X$ and $Y$ as optimally coupled. Any measurable function $\phi^*$ such that $\phi^*(X)\sim\mathbb{Q}$ and $(X,\phi^*(X))$ yields an optimal coupling of $X$ and $Y$, is an optimal transport between $X$ and $Y$. In most cases of interest $\pi^*$ arises as the pushforward measure of $(\mathrm{Id},\phi^*)_{\#}\mathbb{P}$ for some measurable function $\phi^*$. For questions of existence and uniqueness as well as further technical background, see \cite{villani2003} and  \cite{villani2006}.
	
In the case of $W_2$ when an optimal $\phi^*$ exists, there is a simple characterization of the optimality of $\phi^*$; i.e. $\phi^*=\nabla\varphi$ is the gradient of a convex function $\varphi:\real^d\rightarrow\real$.
We will rely on this result. We will also rely on orthonormal transformations of vectors; $\mby=A\mbx$, where $A'A=AA'=I$ and $A$ is a $d\times d$ orthonormal matrix, with $I$ being the identity matrix.  Recall that
$\mby'\mby=\mbx'\mbx$ in such cases.

The paper is based on the idea that if $X$ and $Y$ share the same copula, then $\phi^*$ is a known map. This is a result recently found in \cite{alfonsijourdain2014}. For a multivariate distribution function with associated random variable $X=(X_1,\ldots,X_d)$,
$$\mbox{P}(X_1\leq x_1,\ldots,X_d\leq x_d)=F(x_1,\ldots,x_d)=
C\big(F_1(x_1),\ldots,F_d(x_d)\big),$$
for some probability distribution function $C$ on $[0,1]^d$, and $F_i$ is the marginal distribution function for $X_i$. If it is known that two multivariate probability measures $\mathbb{P}$ and $\mathbb{Q}$ on $\real^d$ share a common copula, then the optimal transport reduces to the $1$--dimensional transport on each dimension. That is, for respective marginal distributions $F_i$ and $G_i$, the optimal transport proceeds along marginal dimensions, so $\phi_i^*(x_i)=G_i^{-1}(F_i(x_i))$ and $\phi^*(\mbx)=(\phi_1^*(x_1),\ldots,\phi_d^*(x_d))$. See \cite{alfonsijourdain2014}.

The authors \cite{alfonsijourdain2014} provide a somewhat technical proof. Here we provide a simple proof, which only requires the invariance property of copulas. We state it here without proof and is to be found in, for example, 
\cite{nelson2006}.

\begin{thm}[Optimal Transport for Common Copulas] Let $(X_1,\ldots,X_d)$ and $(Y_1,\ldots,Y_d)$ be random vectors in $\real^d$ with respective marginal distributions $(F_i)$ and $(G_i)$, for $i=1,\ldots,d$,  and sharing a common copula $C$. The optimal transport is $\phi^*(\mbx)=(\phi_1(x_1),\ldots,\phi_d(x_d))$ and
$\phi_i(x_i)=G_i^{-1}(F_i(x_i))$.
\end{thm}
	
\begin{proof}
		First $G_i^{-1}(F_i(X_i))$ gives the correct marginal distribution for $Y_i$. Second, invoking the invariance property of copulas and the fact that each $\phi_i$ is increasing, $(\phi_1(X_1),\ldots,\phi_d(X_d))$ also has the correct copula $C$; hence, the map produces the correct target distribution. 
		
		Then the fact that the $\phi_i$ are increasing functions of the $x_i$ means the derivative of $\phi^*$, $\Phi^*$, is a nonnegative diagonal matrix, and hence $\phi$ is the derivative of a convex function and so is an optimal map.
\end{proof}

We generalize the result of \cite{alfonsijourdain2014} to the following: Assume $\wX$ and $\wY$ are random variables sharing a common copula. If $X=A\wX$ and $Y=A\wY$ then we will show how to find the optimal map between $X$ and $Y$. The result in \cite{alfonsijourdain2014} is obtained when $A=I$ and hence we provide a generalization of their result. We then use this to find, as an illustration, the optimal map between distributions repesenting a certain class of scale mixture of multivariate normal distributions.

In section 2 we state and prove the main result and illustrate its merits in a number of examples. In section 3 we describe how our work applies to the computation a Wasserstein barycenter problem.


\section{Main Results}
The quadratic Wasserstein metric can be written in terms of an integral of inner products;
	\begin{equation}\label{RCL}
		W^2_2(\mathbb{P},\mathbb{Q})=\int_{\real^{2d}}\left(\mbx'\mbx+\mby'\mby-2\mbx'\mby\right) d\pi(\mathbf{x},\mathbf{y}).
	\end{equation}
	The minimization of the above with respect to $\pi$ is the same as the maximization, with respect to $\pi$, of
\begin{equation}\label{key}
\mbox{E}_\pi X'Y=\int_{\real^{2d}} \mbx'\mby\,d\pi(\mathbf{x},\mathbf{y})=\int_{\real^{2d}} (A\mbx)'(A\mby)\,d\pi(\mathbf{x},\mathbf{y})=\mbox{E}_\pi (AX)'(AY)
\end{equation}
for any orthonormal matrix $A$.	
Hence, If $X$ and $Y$ have the same copula, we will then know the optimal map between $AX$ and $AY$.
Also, vica versa, if $AX$ and $AY$ share a common copula then we know the optimal map between $X$ and $Y$.

For example, suppose $d=2$ and 
let 
\begin{IEEEeqnarray}{RCL}
		A&=&\frac{1}{\sqrt{2}}\left(\begin{matrix}
			1 & 1\\
			1 & -1
		\end{matrix}\right).
	\end{IEEEeqnarray}
Write
$\mbox{E}_\pi X'Y=\mbox{E}_{\pi}(X_1Y_1+X_2Y_2)=\half\mbox{E}_{\pi}[(X_1+X_2)(Y_1+Y_2)+(X_1-X_2)(Y_1-Y_2)];$
i.e. $\mbox{E}_{\pi}\mbx'\mby=\mbox{E}_\pi (A\mbx)' (A\mby).$
Hence, if $X_1+X_2$ is independent of $X_1-X_2$, and $Y_1+Y_2$ is independent of $Y_1-Y_2$, so $AX$ and $AY$ share the same common independent copula, we know the optimal map between $AX$ and $AY$, and hence, by virtue of (\ref{key}), between $X$ and $Y$. Indeed, the optimal map between $AX$ and $AY$ is derived from
$$y_1+y_2=F_{Y_1+Y_2}^{-1}(F_{X_1+X_2}(x_1+x_2))\quad\mbox{and}\quad y_1-y_2=F_{Y_1-Y_2}^{-1}(F_{X_1-X_2}(x_1-x_2)),$$
yielding
$$y_1=\half\left\{F_{Y_1+Y_2}^{-1}(F_{X_1+X_2}(x_1+x_2))+F_{Y_1-Y_2}^{-1}(F_{X_1-X_2}(x_1-x_2))\right\}$$
and
$$y_2=\half\left\{F_{Y_1+Y_2}^{-1}(F_{X_1+X_2}(x_1+x_2))-F_{Y_1-Y_2}^{-1}(F_{X_1-X_2}(x_1-x_2))\right\}.$$
In short, if $\phi_A^*$ is optimal between $AX$ and $AY$, then
$\phi^*(\mbx)=A'\phi^*_A(A\mbx)$
is optimal between $X$ and $Y$.

\paragraph{Example 1} We can use the above to find the optimal map between two bivariate normal distributions, represented by random variables $X$ and $Y$, both with zero means, and with covariance matrices
$$\left(\begin{matrix} 1 & 0 \\ 0 & 1 \end{matrix}\right)\quad\mbox{and}\quad \left(\begin{matrix} 1 & \rho \\ \rho & 1 \end{matrix}\right),$$
respectively. 
For both normal distributions, even for $Y$, $X_1+X_2$ is independent of $X_1-X_2$, and $Y_1+Y_2$ is independent of $Y_1-Y_2$. The marginal distributions for each of these four variables is $\mbox{N}(0,2)$, $\mbox{N}(0,2)$, $\mbox{N}(0,2(1+\rho))$ and $\mbox{N}(0,2(1-\rho))$, respectively. Hence,
the optimal map $\phi^*(\mbx)=(\phi_1^*(\mbx),\phi_2^*(\mbx))$ is given by
$$\phi_1^*(\mbx)=\half\left[x_1\big(\sqrt{1+\rho}+\sqrt{1-\rho}\big)+x_2\big(\sqrt{1+\rho}-\sqrt{1-\rho}\big)\right]$$
and
$$\phi_2^*(\mbx)=\half\left[x_1\big(\sqrt{1+\rho}-\sqrt{1-\rho}\big)+x_2\big(\sqrt{1+\rho}+\sqrt{1-\rho}\big)\right].$$
While a known result, it is remarkably straightforward to show using the Parseval identity.
	

\vspace{0.2in}
\noindent
To write out our main results formally and in generality, let $X\sim\mathbb{P}$ and $Y\sim\mathbb{Q}$ be random variables in $\real^d$, with $\mathbb{P}$ and $\mathbb{Q}$ non--degenerate, with $\phi^*$ the optimal map $\phi^*_{\#}\mathbb{P}=\mathbb{Q}$. Let $A$ be a $d\times{}d$ orthonormal matrix and write $X_A\eqd{}AX$ and $\mathbb{P}_A\triangleq{}A_{\#}\mathbb{P}$, with analogous definitions for $Y_A$ and $\mathbb{Q}_A$. Then we have the following characterization of the optimal transport and the corresponding $W_2$ distance between $\mathbb{P}_A$ and $\mathbb{Q}_A$.
	
\begin{thm}
		For $X_A$ and $Y_A$ defined as above we have
		\begin{itemize}
			\item[(1)] $\phi_A^*(\mbx)=A\,\phi^*(A'\,\mbx)$ is the optimal transport for $X_A\mapsto{}Y_A$ \\
			\item[(2)] $W_2(\mathbb{P}_A,\mathbb{Q}_A)=W_2(\mathbb{P},\mathbb{Q})$
		\end{itemize}
\end{thm}
	
\begin{proof}
		For the first part,  we have
			$A\,\phi^*(A'\,X_A)= A\,\phi^*(A'\,A\,X)= A\,\phi^*(X)=A\,Y\eqd Y_A$,
so we obtain the correct target distribution. 
To ascertain optimality, we need to show the derivative of $\phi_A^*$ is a positive demidefinite matrix, knowing that the derivative of $\phi^*$, written as  $\Phi^*$, is positive semidefinite.
Now
$\partial \phi_A^*/\partial\mbx=A\,\Phi^*(A'\mbx)\,A'$
and it is well known that if $\Phi^*$ is positive semidefinite then $A\Phi^*A'$ is also positive semidefinite.
See, for example,  \cite{hornjohnson2013}.
For the second part, we make use of (\ref{key}). 
	\end{proof}
	
	
	
\paragraph{Example 1 continued} More generally in $\real^d$, suppose $X$ is normal, mean 0 and covariance matrix $\Sigma_1$, and $Y$ is normal with mean 0 and covariance matrix $\Sigma_2$. Suppose also there exists a positive definite symmetric matrix $R$ such that
$R\,\Sigma_1\,R=\Sigma_2.$
Hence, for such $R$, i.e. symmetric and positive definite, there exists an orthonormal $P$ such that
$P\,R\,P'=D,$ where $D$ is positive diagonal.

Hence, taking $\widetilde{X}=PX$ and $\widetilde{Y}=PY$, we see that the covariance matrix of $\widetilde{X}$
is $\widetilde{\Sigma}_1=P\Sigma_1P'$ and the covariance matrix of $\widetilde{Y}$ is
$\widetilde{\Sigma}_2=PR\Sigma_1RP'=D(P\Sigma_1P')D$.
So we see that $\widetilde{X}$ and $\widetilde{Y}$ have a common copula and 
$\widetilde{Y}=D\,\widetilde{X}$
is the optimal transform between $\widetilde{X}$ and $\widetilde{Y}$, and hence $Y=RX$ is optimal; i.e. $\phi^*(\mbx)=R\,\mbx$, between $X$ and $Y$.

The $R$ here is found to be
$R=B(B\Sigma_1B)^{-1/2}B$
where $B=\Sigma_2^{1/2}$. It is easy to show that $R$ is both symmetric and positive definite. While this result is known \cite{olkin1982}; the proof of optimal maps between two normal distributions presented here is remarkably straightforward.

\vspace{0.2in}
\noindent
Here we note such a matrix $R$ is not only an optimal transport between Gaussian distributions, but also between any members of the same location--scale family. Hence this demonstrates that two members of a location--scale family can be seen as common rotations of a common copula. That is, for $X\sim\Sigma_1$ and $Y\sim\Sigma_2$, zero mean members of the same location scale family, we have that $\widetilde{X}=PX$ and $\widetilde{Y}=PY$ share the same copula, for an appropriate rotation $P$. Furthermore, a rotation $P$ satisfies this property for $X$ and $Y$ when $P$ is a matrix of eigenvectors of the optimal transport $R$ between $X$ and $Y$. Then $D$, the matrix of eigenvalues of $R$, is the optimal transport between $\widetilde{X}$ and $\widetilde{Y}$.

More generally we have:
\begin{thm}
If $X$ and $Y$ are in $\real^d$ and there exists a symmetric positive definite $d\times d$ matrix $R$ such that
$R\,\mbox{Cov}X\,R=\mbox{Cov}Y,$
and orthonormal matrix $P$ satisfies $P\,R\,P'=D$, where $D$ is a positive diagonal matrix, then $PX$ and $PY$ share the same copula.
\end{thm}

\paragraph{Example 2} We illustrate with $d=2$ and consider the following family of distributions based on 
$X_1=S_1\,Z_1$ and $X_2=S_2\,Z_2$, with all $(S_j)$ and $(Z_j)$ independent with the $Z_j$'s standard normal and the $S_j$'s positive scalar random variables. So $X_1$ and $X_2$ are independent.
Now consider
$\wX_1=\sqrt{\half}\,(S_1Z_1+S_2Z_2)\quad\mbox{and}\quad \wX_2=\sqrt{\half}(S_1Z_1-S_2Z_2).$
Then $\wX=(\wX_1,\wX_2)$ is a random variable which has a scale mixture of normal representation. That is,
$\wX$ given $S=(S_1,S_2)$, is normal with mean 0 and covariance matrix
$$\Sigma(S_1,S_2)=\half\,\left(\begin{matrix}
S_1^2+S_2^2 &  S_1^2-S_2^2\\ 
S_1^2-S_2^2 & S_1^2+S_2^2
\end{matrix}\right).
$$
Hence,
$f_{\wX}(x_1,x_2)=\int\int \mbox{N}_2\big(\mbx|0,\Sigma(s_1,s_2)\big)\,f_1(s_1)\,f_2(s_2)\,d s_1\,d s_2.$
Another member of this family, $\wY$, would be different with different choices for the $S_1$ and $S_2$, having densities $g_1(s_1)$ and $g_2(s_2)$, respectively. If we denote the family as $\{C(\cdot,\cdot)\}$ then we know the optimal map between two members $C(f_1,f_2)$ and $C(g_1,g_2)$.

	
\section{Application to Barycenters}
	
One application of Wasserstein distances is in the calculation of barycenters, or  ``average'' of a set of distributions, with respect to the Wasserstein distance. The $W_2$--\emph{barycenter problem} is as follows: Given $\mathbb{P}_1,\ldots,\mathbb{P}_n\in\mathcal{P}^2(\real^d)$, square integrable probability distributions on $\real^d$, and weights $(\lambda_i)_{i=1}^n$, such that $\sum_{i=1}^{n}\lambda_i=1$, the barycenter problem is to find  $\bar{\mathbb{Q}}\in\mathcal{P}^2(\real^d)$ such that
\begin{IEEEeqnarray}{RCL}
		\bar{\mathbb{Q}}&\triangleq&\argmin_{\mathbb{Q}\in\mathcal{P}_2}\sum_{i=1}^{n}\lambda_i{}W_2^2(\mathbb{P}_i,\mathbb{Q}).
\end{IEEEeqnarray}
Hence, $\bar{\mathbb{Q}}$ is the $W_2$--barycenter and a minimizer for the barycentric cost functional,
\begin{IEEEeqnarray}{RCL}
		V(\mathbb{Q})&\triangleq&\sum_{i=1}^{n}\lambda_i{}W_2^2(\mathbb{P}_i,\mathbb{Q}).
\end{IEEEeqnarray}
For existence  and uniqueness, see \cite{aguehcarlier2011} and \cite{gouicloubes2017}. 

The authors in \cite{alvarezestebanetal2016}  develop a fixed-point approach to finding barycenters. A functional $G:\mathcal{P}^2(\real^d)\rightarrow\mathcal{P}^2(\real^d)$ is constructed whose fixed points, under basic assumptions, are the barycenters of $(\mathbb{P}_i,\lambda_i)_{i=1}^n $. Let $\mathbb{Q}\in\mathcal{P}^2(\real^d)$. For $i=1,\ldots,n$, let $T_i$ be the optimal map such that $T_i(X)\sim\mathbb{P}_i$ with $X\sim\mathbb{Q}$. Then define
\begin{IEEEeqnarray}{RCL}
		G(\mathbb{Q})&\triangleq&\mathcal{L}\left(\sum_{i=1}^{n}\lambda_i\,T_i(X)\right),
\end{IEEEeqnarray}
where $\mathcal{L}$ denotes the distribution, or ``law'', of the underlying random variable.
\begin{prop}
		(\cite{alvarezestebanetal2016}) If $\mathbb{Q}\in\mathcal{P}_{\mathrm{ac}}^2(\real^d)$ then
			$V(\mathbb{Q})\geq V(G(\mathbb{Q}))+W_2^2(\mathbb{Q},G(\mathbb{Q}))$.
As a consequence, $V(\mathbb{Q})\geq{}V(G(\mathbb{Q}))$, with strict inequality if $\mathbb{Q}\ne{}G(\mathbb{Q})$. For a barycenter $\bar{\mathbb{Q}}$, $G(\bar{\mathbb{Q}})=\bar{\mathbb{Q}}$.
\end{prop}

\noindent	
Recently, \cite{zemelpanaretos2019}, points out that when all the $(\mathbb{P}_i)$ share the same copula, convergence  occurs in one iteration, if the starting distribution  also shares the same copula.
This result follows straightfowardly from \cite{aguehcarlier2011}, and see also \cite{kuang2019}. This result is that
$\mathbb{Q}$ is the barycenter for $(\lambda_i,\mathbb{P}_i)_{i=1}^n$, if and ony if for all $x$ in the support of $\mathbb{Q}$, it is that
$x=\sum_{i=1:n} \lambda_i\,T_i(x)$
where $T_i$ is the optimal map from $\mathbb{Q}$ to $\mathbb{P}_i$.

To see how this works when the $(\mathbb{P}_i)$ share a common copula and also with $\mathbb{Q}$, we note that
$$T_i(x)=\left(\begin{matrix}
P_{i1}^{-1}(Q_1(x_1)) \\
\cdots  \\
P_{id}^{-1}(Q_d(x_d))
\end{matrix}\right)
$$
where $P_{ij}$ is the marginal distribution function for component $j$ from $\mathbb{P}_i$ and likewise for $Q_j$ from $\mathbb{Q}$. Hence, the solution is
$$Q_j^{-1}(u)=\sum_{i=1}^n \lambda_j\,P_{ij}^{-1}(u),$$
combined with the knowledge that $\mathbb{Q}$ has a known copula, completes the description of $\mathbb{Q}$.

If we are now looking for the Barycenter of $(\lambda_i, \mathbb{P}_{i}^A)_{i=1}^n$, where $\mathbb{P}_i^A$ is the measure corresponding to $AX_i$, where $X_i$ has measure $\mathbb{P}_i$ and $A$ is an orthonormal $d\times d$ matrix.

\begin{cor}
Given distributions $(\mathbb{P}_{i}^A)_{i=1}^n$ and $\mathbb{Q}^A$, all with common copulas, the Barycenter $\mathbb{Q}^A$ of $(\lambda_i, \mathbb{P}_{i}^A)_{i=1}^n$ is
$$(Q_j^A)^{-1}(u)=\sum_{i=1}^n \lambda_i\,(P_{ij}^A)^{-1}(u),$$
with $\mathbb{Q}^A$ having the same copula as the $(\mathbb{P}_i^A)$.	The Barycenter for $(\lambda_i,\mathbb{P}_i)$ is then given by $\mathbb{Q}$; where $Y$ and $AY$ have distributions $\mathbb{Q}$ and $\mathbb{Q}^A$, respectively. That is, if $\Psi_\mathbb{Q}(\theta)=\mbox{E}\,e^{-\theta' Y}$ is the generating function for $\mathbb{Q}$, then $\Psi_{\mathbb{Q}^A}(\theta)=\Psi_{\mathbb{Q}}(A'\theta)$.	
\end{cor}

\noindent
Of course, this result works backwards in the sense that we could assume the $(\mathbb{P}_i)$ share a common copula and then look for the Barycenter of $(\mathbb{P}_i^A)$. This is the subject of the next example.

\paragraph{Example 3} Suppose the distributions $(\mathbb{P}_i)$ share the independence copula, but have different marginal distributions. Then the $(\mathbb{P}^A_i)$ will not, in general, share the same copula. 
Assume the $(\mathbb{P}_i)$ are in the same location--scale family and the covariance matrices $(\Sigma_i)$ can be characterized in terms of the eigenvectors of $A$, common to each $\Sigma_i$, and the eigenvalues unique to each $\Sigma_i$. 

The optimal map from distribution $\mathbb{Q}$, with diagonal covariance matrix $D=\mbox{diag}(\sigma_1,\ldots,\sigma_d)$, to $\mathbb{P}_i$ with diagonal covariance matrix $D_i=\mbox{diag}(\sigma_{i1},\ldots,\sigma_{id})$, is given by 
$$T_i(x)=\left(\begin{matrix}
x_1\,\sqrt{\sigma_{i1}/\sigma_1} \\
\cdots  \\
x_d\,\sqrt{\sigma_{id}/\sigma_d}
\end{matrix}\right).$$
Hence, the Barycenter $\mathbb{Q}$ has 
$$\sigma_j=\left(\sum_{i=1:n}\lambda_i\,\sqrt{\sigma_{ij}}\right)^2;\quad\mbox{i.e.}\quad
D=\left(\sum_{i=1}^n\lambda_i\,D_i^{1/2}\right)^2.$$
To now characterize the Barycenter $\mathbb{Q}^A$ for the $(\lambda_i,\mathbb{P}_i^A)$, we note that if $Y$ has distribution $\mathbb{Q}$, then
$\mbox{E}e^{-\theta'Y}=g(\theta'\,D\,\theta)$
for some scalar function $g$. Then $\mbox{E}e^{-\theta' AY}=g(\theta'\,A\,D\,A'\,\theta)$
and hence $\mathbb{Q}^A$ has covariance matrix $ADA'$.

\bibliographystyle{plain}
\bibliography{ArxPart}

\end{document}